\documentclass [oneside]{amsart}
\usepackage {amscd}
\usepackage[dvips,final]{graphics}

\advance\textwidth by +1.0in
\advance\oddsidemargin by -0.5in
\advance\evensidemargin by -0.5in

\def\bc{\begin{center}}
\def\ec{\end{center}}
\def\px{\frac{\partial}{\partial x}}
\def\py{\frac{\partial}{\partial y}}
\def\E{\mathcal{E}}
\def\EC{\mathcal{E}}
\def\D{\mathcal{D}}
\def\M{\mathcal{M}}
\def\L{l}
\def\S{\mathcal{S}}
\def\m{\mbox{m}}

\newcommand{\Z}{{\mathbb Z}}
\newcommand{\C}{{\mathbb C}}
\newcommand{\CP}{{\mathbb {C}P(2)}}
\newcommand{\CE}{{\mathcal C}_{exp}}
\newcommand{\cqd}{\ \hfill\rule[-1mm]{2mm}{2mm}}

\newtheorem{theorem}{Theorem}
\newtheorem{proposition}{Proposition}
\newtheorem{lemma}{Lemma}
\newtheorem{definition}{Definition}
\newtheorem{remark}{Remark}
\newtheorem{corollary}{Corollary}
\newtheorem{example}{Example}

\begin{document}

\title{Multiplicity of Invariant Algebraic Curves and Darboux Integrability}

\author{Jaume Llibre$^1$ and Jorge Vit\'orio Pereira$^2$}
\maketitle
\bc
$^1$ Departament de Matem\`atiques, Universitat Aut\`onoma de Barcelona,\\
08193 -- Bellaterra, Barcelona, Spain. 
email: jllibre@mat.uab.es \\

\bigskip
$^2$ Instituto de Matem\'{a}tica Pura e Aplicada, IMPA, Estrada Dona Castorina, 110 \\
Jardim Bot\^{a}nico, 22460-320 - Rio de Janeiro, RJ, Brasil. 
email : jvp@impa.br 
\ec

\begin{abstract}
We define four different kinds of multiplicity of an invariant algebraic curve for a given polynomial
vector field and investigate their relationships. After taking a closer look at the singularities
and at the line of infinity, we improve the Darboux theory of integrability using these new notions of multiplicity.
\end{abstract}

\section{Introduction}\label{s1}

In this paper we want to contribute toward a better understanding of the fascinating relationships between the integrability of 
a polynomial vector field on the complex plane, and the existence of invariant algebraic solutions for such vector fields.
The seminal work on this subject is Darboux's paper \cite{D} of 1878. He showed how to force integrability from the 
abundance of invariant algebraic curves. More precisely, he proved that a polynomial vector field of degree $d$ with at least 
$d(d+1)/2$ invariant algebraic curves has a first integral. This theory of integrability received also contributions from the 
work of Poincar\'e \cite{Po}, who mainly was interested
in the rational first integrals.  

The subject remained almost forgotten until 1979 when Jouanolou \cite{J} published his Lecture Notes. There he shows that if
the number of invariant algebraic curves of a polynomial system of degree $d$ is at least $[d(d+1)/2] + 2$, then the vector 
field has a rational first integral, and in particular all its solutions are algebraic curves. 

Since then, the subject has been revitalized and studied under different viewpoints. Thus, for instance, the 1983 work of 
Prelle and Singer \cite{PS}, using methods of differential algebra, shows that if a polynomial vector field has an elementary
first integral, then it can be computed using only the invariant algebraic curves; i.e., this kind of first integrals can be 
found using Darboux's approach. Nine years later Singer \cite{S} studied the existence of a wider class of first integrals, 
namely Liouvillian first integrals. Roughly speaking, his main result is that the Liouvillian  first integrals have integrating
factors given by Darbouxian functions.

From a more geometrical point of view the work of Sc\'ardua \cite{B} relates the properties of the holonomy groups of invariant
algebraic curves with the existence of first integrals. Later on Camacho and Sc\'ardua \cite{CS} proposed to find a wider class of
mathematical objects generalizing the Liouvillian functions.

In this paper we adopt a more classical point of view, closer to Darboux's original one, and improve his theory in various aspects.
One of the main ideas coming from the work of Darboux, is the control of the space of cofactors associated to the invariant algebraic curves.
If the polynomial vector field has degree $d$, then the cofactors are polynomials of degree at most $d-1$, hence are contained in the
finite dimensional vector space $\C_{d-1}[x,y]$. 

In \cite{CLS} the authors used singular points to reduce the dimension of the space of possible cofactors and consequently, the number
of conditions to force Darboux integrability. Here we strength this approach considering  not just the singular points, but
also their degree of degeneracy.

In \cite{P}  the extactic curves for polynomials vector fields are introduced. Essentially, they are the curves of inflection
and higher order inflection points for the orbits of the vector field. In the context of integrability they provide new bounds 
for the number of invariant algebraic curves. These bounds take into account the degree of the invariant algebraic curves, and
for a fixed degree they provide, in general, much better bounds than Darboux's classical approach. 

In this work we use the extactic curves to define two notions of algebraic multiplicity. The first (see Subsection \ref{AM})
is more useful to force integrability, while the second (see Section \ref{strong}) turn out to be a more precise tool 
to bound the geometric multiplicity, another notion defined here. Finally we show how the geometric multiplicity gives rise
to exponential factors, introducing the notion of derivative of an invariant algebraic curve. The germ of the idea of the interaction
between exponential factors and geometric multiplicity can be found in Christopher's paper \cite{CH} and in \cite{CL1,CL2}.
Motivated by this interaction we define the notion of integrable multiplicity. This is the key notion which allow us to improve
the Darboux theory of integrability.

The interplay between the different notions of multiplicity and the Darboux integrability can be sumarized in the following diagram.

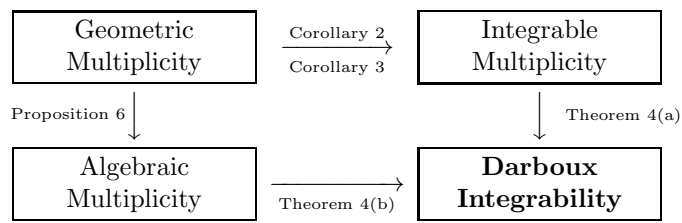
\begin{equation*}
\begin{CD}
\fbox{\begin{minipage}{3cm}\bc Geometric Multiplicity\ec\end{minipage}} 
@>\mbox{{\tiny Corollary \ref{gii}}}>\mbox{{\tiny Corollary \ref{ccl}}}> \fbox{\begin{minipage}{3cm}\bc Integrable Multiplicity\ec\end{minipage}} \\
@V{\mbox{{{\tiny Proposition \ref{ineq1}}}}}VV @VV{\mbox{\tiny \begin{minipage}{2cm}\bc  Theorem \ref{main}(a) \ec\end{minipage}}}V \\
\fbox{\begin{minipage}{3cm}\bc Algebraic Multiplicity\ec\end{minipage}} @>>\mbox{{\tiny Theorem \ref{main}(b)}}> 
\fbox{\begin{minipage}{3cm}\bc{\bf Darboux Integrability}\ec\end{minipage}}
\end{CD}
\end{equation*}

\bigskip

The paper is organized as follows. In Section \ref{s2} we define the basic notions of algebraic invariant curves, 
first integrals, integrating factors, exponential factors and extactic curves. Section \ref{multc} is responsible for introducing
the notions of geometric, algebraic  and integrable multiplicities, besides the concept of derivative of an algebraic invariant
curve. There we also investigate some relations between them. Section \ref{sing} studies the singular points and the line at infinity
preparing the background to Section \ref{darboux}, where we improve the Darboux theory of integrability. Finally in Section 
\ref{strong} we define the strong algebraic multiplicity, and comment how it can be used to provide sharper upper bounds for 
the geometric multiplicity. We conclude with some remarks given in Section \ref{final}.

\bigskip

\tableofcontents

\section{Preliminaries}\label{s2}

\subsection{Vector fields}

\begin{definition}\label{grau}\rm
We say that $X$ is a {\it (polynomial) vector field} of (affine) degree $d$ on $\C^2$ if it can be written in the form
\begin{equation} \label{campo}
  X = a(x,y) \px + b(x,y) \py   \ ,
\end{equation}
where $a(x,y)$ and $b(x,y)$ are polynomials in $\C[x,y]$ such that the maximum degree of $a$ and $b$ is $d$.
\end{definition}

Another definition of the degree of a vector field $X$ is the following one: $X$ has degree $d$ if the number of tangencies
of $X$ with a generic straight line is $d$. In general the two definitions do not coincide, but they differ by at most $1$.
For more details about this see \cite{LS}.

The set of all vector fields of degree $d$ can be seen as a Zariski open set in the affine 
space given by the coefficients of the vector fields.

\subsection{Invariant algebraic curves and exponential factors}

\begin{definition}\rm
Let $f \in \C[x,y]$. The algebraic curve $\{f=0\}$, or simply $f$,  is said to 
be {\it invariant} by the vector field $X$ when $X(f)/f$ is a polynomial.
In that case $L_f = X(f)/f$ is called the {\it cofactor} of the invariant algebraic curve $f$.
\end{definition}

Observe that the degree of $L_f$ is at most the degree of $X$ minus one, and that an invariant algebraic curve of $X$ is formed by solutions 
of the vector field $X$. In other words a solution of $X$ has either empty intersection with the zero set of $f$, or it is entirely contained
in $\{f=0\}$. 

The next two propositions can be found in \cite{CL2}.

\begin{proposition}\label{propfactor}
Let $f \in \C[x,y]$ and $f=f_1^{n_1}\cdots f_r^{n_r}$ be its factorization in irreducible factors. Then, for a vector field $X$, $f$ is
an invariant algebraic curve with cofactor $L_f$ if, and only if $f_i$ is an invariant algebraic curve for each $i=1,\ldots,r$ with cofactor
$L_{f_i}$. Moreover $L_f = n_1L_{f_1} + \cdots + n_rL_{f_r}$.
\end{proposition}

\begin{definition} \label{exp} \rm
Let $f, g \in \C[x,y]$, we say that $e=\exp(g/f)$ is an {\it exponential factor} of the vector field $X$ of degree $d$, if $X(e)/e$
is a polynomial of degree at most $d-1$. This polynomial is called the {\it cofactor} of the exponential factor $e$, which
we denote by $L_e$. The quotient $g/f$ is an {\it exponential coefficient} of $X$.
\end{definition}

Note that despite the fact that an exponential factor satisfies the same kind of equation than an invariant algebraic curve, a priori, 
it has no relation with the solutions of the vector field because its zero set is empty. Although we have the following result.

\begin{proposition}
If $e = \exp\left(g/f\right)$ is an exponential factor for the vector field $X$, then $f$
 is an invariant algebraic curve and $g$ satisfies the equation
 \[
   X(g) = gL_f + fL_e \ .
 \]
\end{proposition}

\subsection{First Integrals and Integrating Factors}

This work is mainly interested in two classes of first integrals, the rational ones, and the {\it Darbouxian}; that is a product 
of complex powers of polynomials and complex powers of exponential of rational functions. We recall that by first integral we mean:

\begin{definition} \rm
Let $U$ be an open subset of $\C^2$. We say that a non constant (multi--valued) function $H : U \to \C$ is a {\it first integral}
of a vector field $X$ on $U$ if, and only if, $X|_U(H) = 0$. When $H$ is the restriction of a rational, respectively Darbouxian, function 
to $U$ then we say that $H$ is a rational, respectively Darbouxian, first integral.
\end{definition}

Besides the first integrals we are also interested in the integrating factors. Which can be formally defined as follows.

\begin{definition} \rm
We say that a non constant (multi--valued) function $R : U \to \C$ is an {\it integrating factor} of a vector field $X$
on $U$ if, and only if, $X|_U(R) = - \mbox{div} X|_U \cdot R$. We recall that $\mbox{div}$ denotes the divergence of the vector field.
\end{definition}

If we know an integrating factor we can compute a first integral. Reciprocally if $H$ is a first integral of the vector field 
(\ref{campo}), then there is a unique integrating factor $R$ satisfying 
\begin{equation}\label{intfactor}
  R a = \frac{\partial H}{\partial y}  \quad \mbox{and} \quad R b = -\frac{\partial H}{\partial x} \ .
\end{equation}
Such $R$ is called the {\it integrating factor associated to $H$}.

\subsection{Extactic curves}

The definitions and results of this subsection comes from the work of the second author. Here we state and prove simplified versions
adapted to the complex plane.

Intuitively the objects defined below measure the order of contact of the solutions 
of a vector field with algebraic curves of given degree. For more information see \cite{P}.

\begin{definition}\rm
If $X$ is a vector field on $\C^2$ the {\it $n$--th extactic curve of $X$}, $\E_n(X)$, 
is given by the equation   
\begin{equation}\label{extmat}
\mbox{det} \begin{pmatrix}       v_1    &    v_2       & \cdots & v_l          \cr
                           X(v_1)       & X(v_2)       & \cdots & X(v_l)       \cr
                           \vdots       & \vdots       & \cdots & \vdots       \cr
                           X^{l-1}(v_1) & X^{l-1}(v_2) & \cdots & X^{l-1}(v_l) \end{pmatrix},
\end{equation}
where  $v_1, v_2, \ldots, v_l $ 
is a basis of the $\C$--vector space of polynomials in $\C[x,y]$ of degree at most $n$ (hence $l=(n+1)(n+2)/2$), 
$X^0(v_i)= v_i$ and $X^j(v_i)= X^{j-1}(X(v_i))$.
\end{definition}

Observe that the definition of extactic curve is independent of the
chosen basis of the $\C$--vector space of polynomials of degree at most $n$.

\begin{theorem} \label{ddd}
Let $X$ be a vector field on $\C^2$. Then $\EC_n(X) = 0$ and $\EC_{n-1}(X) \neq 0$ if, and only if, 
$X$ admits a rational first integral of exact degree $n$.
\end{theorem}

\noindent{\it Proof}:
Let $p \in \C^2$ be a non--singular point of $X$. We may assume 
that $p$ is the origin of $\C^2$.  Suppose that the solution passing through it is parametrized, locally, by $(x,y(x))$. 
Since $\EC_n(X)$ vanishes identically, the composition of our local solution with 
the $n$--Veronese map,
\[
  (x,y) \to (x^n,x^{n-1}y,x^{n-2}y^2,\ldots,y^n,x^{n-1},\ldots,y^{n-1},\ldots,x^2,xy,y^2,x,y,1)
\]
is contained in a hyperplane, therefore $(x,y(x))$ must be contained in an 
algebraic curve of degree less than or equal to $n$. The fact that $\EC_{n-1}(X) \neq 0$ implies  
that the generic solution is of degree at least $n$. This completes the ``only if'' part of the 
statement.

If $X$ admits a first integral of degree $n$ then every invariant curve is of degree at most $n$ 
and hence every point is a $n$--inflection point, i.e., $\EC_n(X) = 0$. Since not every 
invariant curve has degree $n-1$,  $\EC_{n-1}(X) \neq  0$. 
\cqd

\begin{proposition}\label{pext}
Every algebraic curve of degree $n$  invariant by the vector field $X$ is a factor of $\EC_n(X)$.
\end{proposition}

\noindent{\it Proof}: Let $f$ be an invariant algebraic curve of degree $n$.
As it was observed the choice of the basis of the $\C$--vector space plays no role in the definitions
of extactic curve, therefore we can take $v_1=f$. Since 
\[
\begin{array}{lclcl}
X(f)   &=& L_f f \ ,     & &                                     \\
X^2(f) &=& X(L_f f)      &=& \left( L_f^2 + X(L_f) \right) f \ , \\
X^k(f) &=& X(X^{k-1}(f)) &=& \left ( \, \mbox{polynomial} \, \right) f         \ ,
\end{array}
\]
$f$ is a factor of $\EC_n(X)$. 
\cqd

\begin{example}\label{vulpe1} \rm
If we consider the $2$--parameter of vector fields 
\[
  X_{(t,b)} = \left( 2t^2 - 2x^2 \right) \px + \left( b - 4xy - 2t^3y^2 \right) \py \ ,
\]
we have that the first extactic curve of $X_{(t,b)}$ is described by 
\begin{eqnarray*}
  & & y\left( x-t \right) \left( x+t \right) \left( x+t^3y - t\sqrt{\frac{2+bt}{2}} \right) 
  \left( x+t^3y + t\sqrt{\frac{2+bt}{2}} \right) =
  y\L_t^{(1)}\L_t^{(2)}\L_t^{(3)}\L_t^{(4)}.
\end{eqnarray*}
It follows from Proposition \ref{pext} that every invariant straight line, must be contained
in the first extactic curve. It can be easily verified that $y$ is not invariant by any
vector field of the form $X_{(t,b)}$ with $b \ne 0$, and that for every $i \in \{1,2,3,4\}$, $\L_t^{(i)}$ is an invariant
straight line for $X_{(t,b)}$.

The above family of vector field appears in the paper \cite{SV}. \qed
\end{example}

\begin{example}\label{dana1} \rm
The $3$--parameter family of vector fields given by
\[
  X_{(t,b,d)} = \left( -y -bx^2 -dy^2 \right) \px + \left( x + (t-b) xy \right) \py \ ,
\]
corresponds to a quadratic polynomial vector field having a center at origin. It has been studied
in \cite{Sc2}.

The first and the second extactic curves, $\E_1(X_{(t,b,d)})$ and $\E_2(X_{(t,b,d)})$, for $X_{(t,b,d)}$ are
\begin{eqnarray}
& & (1+(t-b)y)\cdot(tbx^4-d^2y^4-2dy^3+(-td+bd)x^2y^2-2dx^2y-x^2-y^2) \ , \label{100} \\
& & (1+(t-b)y)^4\cdot x \cdot (d-t+b(2t-2d)y+tb(b+t)x^2+bd(t+b)y^2) \cdot p_{16}  \  , \label{101} 
\end{eqnarray}
respectively. Here $p_{16}$ denotes a polynomial of degree $16$. One can see that  the linear factor of (\ref{100}) and the quadratic factor
of (\ref{101}) are the unique invariant algebraic curves of the respective degree for the vector field $X_{(t,b,d)}$.
\qed
\end{example}

\section{Notions of multiplicity for invariant algebraic curves}\label{multc}

\subsection{Geometric Multiplicity}\label{GM}

Let $X$ be a vector field on $\C^2$ of degree $d$.

\begin{definition}\label{multgeo}\rm
An invariant algebraic curve $f$ of degree $n$ for the vector field $X$, has {\it $l$--geometric multiplicity $m$ } if,
there exists a sequence of vector fields $X_k$, of degree $d$, converging to $X$, such that each $X_k$
has $m$ distinct invariant algebraic curves of degree at most $l$, $f_k^{(1)}, \ldots , f_k^{(m)}$, converging to $f$ as $k$ goes 
to infinity. Here, by distinct we mean that pairwise they have no common factors. We denote $m$ by $\mu_{g,l}(X,f)$, or if no confusion can arise $\mu_{g,l}(f)$. 

Finally we define the {\it geometric multiplicity}, $\mu_{g}(X,f)$, of the invariant algebraic curve $f$ with respect to $X$ as 
\[
  \lim_{l \to \infty} \mu_{g,l}(X,f) \ .
\]
\end{definition}

\begin{remark}\rm
By convergence  we mean the convergence of the coefficients of the polynomials involved,
either in the case of vector fields, or in the case of algebraic invariant curves.
\end{remark}

\begin{remark}\rm
In Definition \ref{multgeo} if $l < n$, then clearly $\mu_{g,l}(X,f) = 0$.
\end{remark}

Note that in the definition of geometric multiplicity it is allowed that it becomes infinity. 
But if this is the case we will show  in Proposition \ref{inf}  
that the vector field has a rational first integral.

\begin{example}\label{vulpe1geo} \rm
To illustrate the geometric multiplicity of an invariant algebraic curve we shall compute the multiplicity
of the invariant straight line $\L = \{ x=0 \}$  for the vector field $X_{(0,b)}$ defined in Example \ref{vulpe1}.
From that example we know that for any $t \neq 0$, $X_{(t,b)}$ has exactly four distinct invariant straight lines, and hence one can
conclude that $\mu_{g,1}(X_{(0,b)},L) \ge 4$. Later on we will see that $\mu_{g,1}(X_{(0,b)},L) = 4$. 
\qed
\end{example}

\begin{example}\label{dana1geo} \rm
For the vector fields $X_{(0,b,d)}$ given in Example \ref{dana1}, the invariant straight line $1 - by$ raised to 
the second power has $2$--geometric multiplicity at least $2$. This is due to the fact that if we fix $b$ and $d$ and make $t$ tend to zero,
it is possible to check that we have the invariant conic $C = d-t+b(2t-2d)y+tb(b+t)x^2+bd(t+b)y^2$ tends to $(1-by)^2$. For the moment we cannot
give upper bounds to the $l$--geometric multiplicity, that, we will do in the next subsection. 
\qed
\end{example}

\begin{proposition}\label{familia}
Let $f$ be an algebraic curve of degree $n$ invariant by the vector field $X$. 
Whenever $\mbox{$\mu_{g,l}(X,f) = m$}$, there exists an analytic $1$--parameter family of polynomial vector 
fields with degree $d$, $X_t$ where $\mbox{$t \in (-\epsilon,\epsilon)$}$, 
analytic in $t$ such that for each $\mbox{$t \in (-\epsilon,\epsilon) \setminus \{0\}$}$, $X_t$
has $m$ distinct invariant algebraic curves $f_t^{(1)},\ldots, f_t^{(m)}$ of degree at most $l$
converging to $f$ as $t$ tends to zero, and $X_0=X$.
\end{proposition}

\noindent{\it Proof}: It follows easily from the fact that the set of vector fields with invariant algebraic curves of degree
at most $n$ is a Zariski closed set in the space of vector fields of degree $d$.
\cqd

\begin{proposition}\label{inf}
Let $f$ be an invariant algebraic curve for the vector field $X$. If $\mu_{g}(X,f) = \infty$, then $X$ admits a rational first integral.
\end{proposition}

\noindent{\it Proof}: Suppose $X$ has degree $d$. Since $\mu_{g}(X,f) = \infty$, there exists an $l$ such that
\[
\mu_{g,l}(X,f) > \eta = \frac{d(d+1)}{2} + 1 \ . 
\]
From Proposition \ref{familia} there exists an analytic $1$--parameter family $X_t$,
$t \in (-\epsilon,\epsilon)$, such that for $t \neq 0$, $X_t$ has at least $\eta$ distinct invariant algebraic curves.
It is known from \cite{J} (see a very short proof in \cite{CL2}) that a polynomial vector field of degree $d$ 
having more than $\eta$ distinct invariant algebraic curves has a rational first integral. Hence if $t \neq 0$, $X_t$ admits a rational first integral. 

Let $\mathcal{S}_i$ be the Zariski closed set in the space of all vector fields of degree $d$ with rational first integral of degree
at most $i$. The intersection of $\mathcal{S}_i$ with the analytic family $X_t$ is either all the family or is a discrete set. In the first case 
the proposition is proved. So suppose that for all $i$, the intersection of $\mathcal{S}_i$ with the family $X_t$ is a discrete set, and taking
their union we obtain an enumerable set. 
This is in contradiction with the fact this set must be equal to $X_t$, with $t \in (-\epsilon,\epsilon) \setminus \{0\}$.
\cqd

\vskip 0.2cm

We point out that Schlomiuk in \cite{Sc} defines a notion of geometric multiplicity of an invariant algebraic curve with respect to 
a family of polynomial vector fields. There, besides the use of the topology of the space of coeffiecients, as we do, she uses a kind
of  Hausdorff topology, see \cite{Sc}. We also like to remark that her definition is with respect to a fixed family,
while ours do not depend on the choice of any particular family. These turn out to be the essential differences between the two notions. 
They are both relevant, and we choose the present approach 
because together with the notion of algebraic multiplicity we can provide computable upper bounds of it without considering
any particular family.

To make more clear the difference observe that in Example \ref{dana1geo}, for Schlomiuk the invariant straight line has
multiplicity $2$, with respect to the considered family  .
For us the invariant straight line has $2$--geometric multiplicity $1$ (see Example \ref{dana1alg}),
and its square has $2$--geometric multiplicity exactly $2$, see Example \ref{decide}.

\subsection{Algebraic Multiplicity}\label{AM}

\begin{definition}\rm
An invariant algebraic curve $f$ of degree $n$ for the vector field $X$ has {\it $l$--algebraic multiplicity $m$} when 
$m$ is the greatest positive integer such that the
$m$--th power of $f$ divides $\E_l(X)$. We denote $m$ by $\mu_{a,l}(X,f)$, or if the context is clear $\mu_{a,l}(f)$.
\end{definition}

As in the definition of $l$--geometric multiplicity, if $l < n$ then $\mu_{a,l}(X,f)=0$.

\begin{proposition}\label{ineq1}
The following inequality holds :
\[
  \mu_{g,l}(X,f) \leq \mu_{a,l}(X,f) \ .
\]
\end{proposition}

\noindent{\it Proof}: Suppose that the degree of $f$ is $n$, with $\mu_{g,l}(X,f)=m$. 
Let $X_k$ be a sequence of vector fields
of degree $d$ converging to $X_0=X$, with $m$ distinct invariant algebraic curves $f_k^{(1)}, \ldots, f_k^{(m)}$ converging to $f$.
By Proposition \ref{pext} the product, $f_k^{(1)} \cdots f_k^{(m)}$
divides $\E_l(X_k)$. Since  $f_k^{(1)} \cdots f_k^{(m)}$ are distinct,
\[
 \lim_{k \to \infty}{f_k^{(1)} \cdots f_k^{(m)}} = f^m \ ,
\]
implies that $f^m$ divides $\E_l(X)$. 
\cqd

\begin{example}\label{vulpe1alg} \rm
In Example \ref{vulpe1} we have calculated the extactic curve of $X_{(t,b)}$. From those calculations it follows
\begin{eqnarray*}
  \mu_{a,1}(X_{(t,b)},\L_t^{(i)}) &=&  1 , \mbox{ for } t\neq 0 \mbox{ and } i = 1, 2, 3, 4 \ ; \\
  \mu_{a,1}(X_{(0,b)},\L      ) &=&  4 \ .
\end{eqnarray*}
Therefore, by Proposition \ref{ineq1}, we obtain the equality 
\[
   \mu_{g,1}(X_{(0,b)},\L) = \mu_{a,1}(X_{(0,b)},\L) = 4
\]
mentioned in Example \ref{vulpe1geo}.
For $n \ge 2$ we have  $\mu_{a,n}(X_{(0,b)},\L^n) \ge 4$ by Proposition \ref{ineq1} and Example \ref{vulpe1geo}, although it can be much more.
Using a computer algebra system one can verify that 
\begin{eqnarray*}
  \E_2(X_{(t,b)})&=&by(x-t)^4(x+t)^4(-2t^2+2x^2-bt^3+4t^3xy+2t^6y^2)^4\cdot \\
                 & &(90t^3xy^3-18bt^3y^2+18t^2y^2+90x^2y^2-45bxy+5b^2) \ ,
\end{eqnarray*}
consequently $\mu_{a,2}(\L^2) = 8$. \qed
\end{example}

\begin{example}\label{dana1alg} \rm
Now we study the algebraic multiplicity of the invariant straight line for the vector field $X_{(0,b,d)}$ given in Example \ref{dana1} and 
analyzed in Example \ref{dana1geo} under the prism of geometric multiplicity. 
From the expression of the first extactic curve $\E_1(X_{(t,b,d)})$ we get immediately that $\mu_{a,1}(X_{(0,b,d)},1-by) = 1$, and consequently
by Proposition \ref{ineq1}, $\mu_{g,1}(X_{(0,b,d)},1-by) = 1$.

Similarly using the expression of the second extactic curve $\E_2(X_{(t,b,d)})$ it is easy to conclude that $\mu_{a,2}(X_{(0,b,d)},(1-by)^2) = 3$.
From Example \ref{dana1geo} we obtain
\[
  2 \le  \mu_{g,2}(X_{(0,b,d)},(1-by)^2) \le 3 \ .
\]
\qed
\end{example}

Observe that in \cite{Sc}, Schlomiuk also define a notion of algebraic multiplicity. In this case this notion is
of complete different nature than ours. Although, there she shows that his notions of geometric and algebraic multiplicity coincides
in the particular case of quadratic Hamiltonian vector fields. 

\subsection{Integrable Multiplicity}\label{IM}

\subsubsection{Derivative of an invariant algebraic curve}

\begin{definition} \label{der}\rm
Let $X_t$, $t \in (-\epsilon,\epsilon)$, be an analytic $1$--parameter family of vector fields of degree $d$ such that, for each 
$t \neq 0$, $X_t$ has two distinct invariant algebraic curves $f_t^{(1)}$ and $f_t^{(2)}$. Let $f$ be an invariant algebraic curve of
$X = X_0$. Assume that $f_t^{(1)},f_t^{(2)} \to f$ as $t \to 0$, and that 
\begin{eqnarray}
  \lim _{t \to 0} \frac{f_t^{(2)} - f_t^{(1)}}{t^j} &=& 0 \quad \mbox{ for } j = 1, \ldots, r-1 \ ,\label{derivada1} \ \\
  \lim _{t \to 0} \frac{f_t^{(2)} - f_t^{(1)}}{t^r} &\neq& 0  \ . \label{derivada2}
\end{eqnarray}
In that case we call the limit (\ref{derivada2}) a {\it derivative} of the invariant algebraic curve $f$. From now on we denote the 
$\C$--vector space generated by all derivatives of $f$ by $\mathcal{D} f$.
\end{definition}

Although the next series of results are new, the germs of its demonstrations are already presented in Section 5 of \cite{CL2}.

\begin{lemma}
Using the notation of Definition \ref{der}, if (\ref{derivada1}) and (\ref{derivada2}) hold we have that
\[
  L_t^{(2)} = L_t^{(1)} + t^rL + O(t^{n+1}) \ ,
\]
where $L_t^{(i)}$, $i=1,2$, are the cofactors of $f_t^{(i)}$ with respect to the vector field $X_t$, and $L$ is a polynomial of degree at most $d-1$.
\end{lemma}

\noindent {\it Proof}: From the following computation
\[
X \left(  \frac{f_t^{(2)} - f_t^{(1)}}{t^j} \right) =  \frac{L_t^{(2)}f_t^{(2)} - L_t^{(1)}f_t^{(1)}}{t^j} 
= \frac{ \left( L_t^{(2)} - L_t^{(1)} \right) f_t^{(1)}}{t^j} + O(t^{r-j}) \ ,
\]
and (\ref{derivada1}) we get 
\[
 \lim_{t \to 0} \frac{  L_t^{(2)} - L_t^{(1)} }{t^j} f = 0 \ ,
\]
whenever $j < r$. Therefore, we can write 
\[
 L_t^{(2)} = L_t^{(1)} + t^rL + O(t^{r+1}) \ ,
\]
and since  $L_t^{(i)}$, $i=1,2$, are polynomials of degree at most $d-1$ it follows that $L$ has degree at most $d-1$.
\cqd
  
\begin{theorem}\label{conta}
Let $X$ be a vector field of degree $d$. If $g$ is a derivative of an invariant algebraic curve $f$ of $X$, then $g/f$ is an exponential coefficient  
of $X$.
\end{theorem}

\noindent{\it Proof}: Since $g$ is a derivative of $f$ there exists an analytic $1$--parameter 
family of vector fields $X_t$ of degree $d$, 
such that $X_0 = X$ and for a convenient $\epsilon$ there exists for 
every $t \in (-\epsilon,\epsilon)$  two distinct algebraic curves
$f_t^{(1)}$, $f_t^{(2)}$ invariant by $X_t$, both converging to $f$, and satisfying 
\[
   \lim _{t \to 0} \frac{f_t^{(1)} - f_t^{(2)}}{t} = g \ .
\]
We will need the next formulas later on:
\begin{eqnarray}
 X_t \left( \frac{f_t^{(2)}}{f_t^{(1)}} \right) &=&  
  \frac{X_t\left(f_t^{(2)}\right)f_t^{(1)} - X_t\left(f_t^{(1)}\right)f_t^{(2)}}{\left(f_t^{(1)}\right)^2} \nonumber \\
  &=& \left( L_t^{(2)} - L_t^{(1)} \right) \frac{f_t^{(2)}}{f_t^{(1)}} \nonumber \\
  &=& t^rL\frac{f_t^{(2)}}{f_t^{(1)}} + O(t^{r+1}) \label{formula} \ .
\end{eqnarray}                                                      

We define $g_t$ as
\[
  \frac{f_t^{(2)} - f_t^{(1)}}{t^r} \ .
\]
From (\ref{formula}) it follows that 
\begin{eqnarray*}
X_t \left( \left( \frac{ f_t^{(2)} }{ f_t^{(1)} } \right)^{1/t^r } \right) &=&  
\frac{1}{t^r}\left( \frac{ f_t^{(2)} }{ f_t^{(1)} } \right)^{1/t^r}
\left( \frac{f_t^{(2)}}{f_t^{(1)}} \right) ^{-1}X_t\left(\frac{f_t^{(2)}}{f_t^{(1)}} \right) \\
 &=&  \left( \frac{f_t^{(1)}+t^{r}g_t}{f_t^{(1)}} \right)^{1/t^r}
\left( \frac{f_t^{(1)}+t^{r}g_t}{f_t^{(1)}} \right)^{-1}
\left( L\frac{f_t^{(2)}}{f_t^{(1)}} + O\left(t\right)\right) \\
 &=& \left( 1 + t^r\frac{g_t}{f_t^{(1)}}\right)^{1/t^r}\big( 1 + O\left(t^r\right) \big)
\left( L\frac{f_t^{(2)}}{f_t^{(1)}} + O\left(t\right)\right) \ .
\end{eqnarray*}
These equalities imply that
\[
X_t \left( \left( 1 + t^r \frac{g_t}{f_t^{(1)} } \right)^{1/t^r}\right) = 
\left( L\frac{f_t^{(2)}}{f_t^{(1)}} + O\left(t\right)\right)\left( 1 + t^r \frac{g_t}{f_t^{(1)} } \right)^{1/t^r} \ .
\]
Finally taking  the limit in this last equality when $t$ tends to zero we obtain
\[
X\left( \exp (g/f) \right) = L\exp (g/f) \ .
\]
Consequently $g/f$ is an exponential coefficient of $X$, and this completes the proof of the theorem.
\cqd

\begin{corollary}\label{xxx}
Given an invariant algebraic curve $f$ of the vector field $X$ with $\mu_g(X,f) > 1$, 
then exists at least one exponential coefficient having non--trivial 
denominator $f_1$, such that $f_1$ divides $f$.
\end{corollary}

To clarify the notions of derivative of an invariant algebraic curve we now study some examples.

\begin{example} \rm
Consider the $2$--parameter family $X_{(t,b)}$ given in Example \ref{vulpe1}. In Example \ref{vulpe1geo} we showed that $\mu_{g,1}(X_{(0,b)},x)=4$. 
The vector field $X_{(t,b)}$ has exactly four invariant straight lines, which we called $\L_t^{(i)}$ for $i=1,\ldots,4$. To obtain derivatives of the 
line $x$ we have to compute the limits 
\[
  \lim_{t \to 0} \frac{\L_t^{(i)} - \L_t^{(j)}}{t} \ , \quad \mbox{for} \, \, i, j \in \{1,2,3,4 \} \quad \mbox{and} \, \, i\neq j \ .
\]
Doing such computations we have
\begin{eqnarray}
   \lim_{t \to 0} \frac{\L_t^{(2)} - \L_t^{(1)}}{t} &=& \lim_{t \to 0} \frac{ (x+t) - (x-t)}{t} = 2 \ , \label{eq21} \\
   \lim_{t \to 0} \frac{\L_t^{(3)} - \L_t^{(1)}}{t} &=& \lim_{t \to 0} \frac{\left( x+t^3y - t\sqrt{1 + bt/2} \right)  - (x-t)}{t} = 0 \ , \label{eq31}\\
   \lim_{t \to 0} \frac{\L_t^{(4)} - \L_t^{(1)}}{t} &=& \lim_{t \to 0} \frac{\left( x+t^3y + t\sqrt{1 + bt/2} \right)  - (x-t)}{t} = 2 \ . \label{eq41}  
\end{eqnarray}
From equation (\ref{eq21}) or (\ref{eq41}) we obtain the same exponential coefficient, namely $2/x$. Since the limit (\ref{eq31}) is zero, we compute
\[
  \lim_{t \to 0} \frac{\L_t^{(3)} - \L_t^{(1)}}{t^2} = \lim_{t \to 0} \frac{\left( x+t^3y - t\sqrt{1 + bt/2} \right)  - (x-t)}{t^2} = -\frac{b}{4} \ .
\]
And from that we obtain the exponential coefficient $-b/(4x)$. We note that this last exponential
coefficient only differs by multiplication by a constant from the previous one, and hence they provide essentially a unique exponential factor.

To find others exponential coefficients we consider products of two invariant straight lines and do similar computations. For example,
if we take the products
of $\L_t^{(1)}\L_t^{(3)}$ and $\L_t^{(2)}\L_t^{(4)}$, and calculate 
\[
  \lim_{t \to 0} \frac{\L_t^{(1)}\L_t^{(3)} - \L_t^{(2)}\L_t^{(4)}}{t} = -4x \ ,
\]
we obtain as exponential coefficient $ -4x / x^2 = -4/x $, which again gives the same exponential factor as before. On the other hand if we compute
\[
\lim_{t \to 0} \frac{\L_t^{(1)}\L_t^{(2)} - \L_t^{(3)}\L_t^{(4)}}{t^3} = -2xy + \frac{b}{2}  \ ,
\]
and we get as exponential coefficient $(-2xy + b/2)/x^2$, which gives rise to a new exponential factor. \qed
\end{example}

\begin{example} \rm
Now, we want to study the derivative of the square of the invariant straight line for the vector field $X_{(0,b,d)}$, already studied in
Examples \ref{dana1}, \ref{dana1geo} and \ref{dana1alg}. The careful reader must have noted that in fact the conic do not tend to $(1-by)^2$, 
but in fact tend to $d \cdot (1-by)^2$ when $t$ tends to zero. So from the definition of derivative of an multiple invariant curve, we obtain that
one of the derivatives (this can be the only one) of $(1-by)^2$ is given by
\[
  g = \lim_{t \to 0} \frac{d(1-by)^2 - C_{(t,b,d)}}{t} = 1  + 2(d -b)y   - b^2x^2 - 3bdy^2 \ .
\]
Therefore $g/(1-by)^2$ is an exponential coefficient of $X_{(0,b,d)}$.
\qed 
\end{example}

\subsubsection{Notion of integrable multiplicity}

\begin{definition}\label{coefexp} \rm
Let $f$ be an invariant algebraic curve of the vector field $X$. We denote by $\CE(f)$, the $\C$--vector subspace of the 
$\C$--vector space of rational functions on $\C^2$, formed by all exponential coefficients of $X$ of the form $g/f^k$ where
$g \in \C[x,y]$ and $k$ runs over all the non--negative integers.
\end{definition}

\begin{proposition}
The following inclusion holds
\[
  \CE(f) \supseteq \left\langle \frac{\D f}{f},\frac{\D f^2}{f^2},\ldots,\frac{\D f^n}{f^n},\ldots \right\rangle \ .
\]  
Here the right hand side of the previous expression denotes the minimal $\C$--vector subspace of the space of rational 
functions that contains all the vector spaces $\D f^n/f^n$ for any positive integer. 
\end{proposition}

\noindent{\it Proof}: It is an easy consequence of Theorem \ref{conta}. 
\cqd

\begin{definition}\label{multint}\rm
An invariant algebraic curve $f$ for the vector field $X$, has {\it integrable multiplicity $m$ } if
\[
 \dim_{\C}\CE(f) = m - 1 \ .
\]
We denote $m$ by $\mu_i(X,f)$, or if no confusion is possible $\mu_i(f)$. 
\end{definition}

The minus one in the above equality is due to the fact that we also take into account the invariant algebraic curve $f$ in the
notion of integrable multiplicity.

With this new notions, we are able to reformulate Corollary \ref{xxx}, as follows.

\begin{corollary}\label{gii}
Let $f$ be an invariant algebraic curve of the vector field $X$. If $\mu_g(f) > 1$, then $\mu_i(f) > 1$.
\end{corollary}

\subsection{Relations between the different notions of multiplicity}

The next result can be find as Theorem 6.3 in \cite{CL1}. We include it here, because it shows how to force geometric multiplicity
from integrable multiplicity.

\begin{proposition}\label{cl}
Suppose a vector field $X$ has an exponential factor $ \exp(g/f)$, where $g=0$ and $f=0$ are non-singular and have normal
crossings with themselves and the line at infinity. Then the vector field is of the form
\begin{equation}\label{ff}
  \left( a_of^2 - a_1ff_y - a_2(g_yf - gf_y) \right)\px + \left( a_3f^2 + a_1ff_x + a_2(g_xf - gf_x) \right) \py \ , 
\end{equation}
where $\deg(a_i)$ satisfy the expected bounds. In particular, the exponential can be seen to be the limit of the invariant curves
$f$ and $f + \epsilon g$ as $\epsilon$ tends to zero in the family of vector fields $X_\epsilon=p_\epsilon \px + q_\epsilon \py$, where
\begin{eqnarray*}
  p_\epsilon &=& a_0f(f + \epsilon g) - (a_1 + \epsilon^{-1}a_2)f(f + \epsilon g)_y + \epsilon^{-1}a_2f_y(f + \epsilon g) \ , \\
  q_\epsilon &=& a_3f(f + \epsilon g) + (a_1 + \epsilon^{-1}a_2)f(f + \epsilon g)_x - \epsilon^{-1}a_2f_x(f+ \epsilon g)  \ ,
\end{eqnarray*}
which tends to (\ref{ff}) as $\epsilon$ tends to zero.
\end{proposition}

From Proposition \ref{cl} we get easily the next corollary.

\begin{corollary}\label{ccl}
If $f$ is an invariant algebraic curve of a vector field $X$ and $g$ is a derivative of $f$ such that $f$ and $g$ satisfies the hypothesis
of Proposition \ref{cl},  then $\mu_g(X,f) \ge 2$.
\end{corollary}

We do not know if it is true or not that $\mu_i(X,f) = \sup_n \mu_g(X,f^n)$.

\section{The role played by singular points and the line at infinity} \label{sing}

\subsection{Finite singular points}

Let $I \subset \C[x,y]$ be an ideal and $p = (x_0,y_0)$ a point of $\C^2$. Denote by $TI_p$ the ideal generated by all monomials
$(x-x_0)^i(y-y_0)^j$ which appears with nonzero coefficient in the Taylor series expansion around $p$ of some polynomial $f \in I$.

Now, let $X$ be a vector field of degree $d$, written in the form (\ref{campo}), and $I$ be ideal generated by $a$ and $b$.
Consider the following exact sequence of $\C$--vector spaces
\begin{equation*}
\begin{CD}
  0 @>>> \ker \pi_p @>>> \C_{d-1}[x,y] @>\pi_p>> \displaystyle{\frac{\C_{d-1}[x,y]}{TI_p \cap \C_{d-1}[x,y]}}@>>> 0 \ ,
\end{CD}
\end{equation*}
where $\pi_p$ is the natural projection and $\C_{d-1}[x,y]$ is the $\C$--vector space of polynomials of degree at most $d-1$.
Note that $\ker \pi_p = \C_{d-1}[x,y]$ if, and only if, $p$ is a non-singular point of $X$.

\begin{proposition} \label{aaa}
If $f$ is an invariant algebraic curve of $X$ not passing through $p$, then $L_f \in \ker \pi_p$.
\end{proposition}

\noindent{\it Proof}: Without loss of generality we can assume that $p=(0,0)$.
Since $f$ do not pass through $p$, after multiplying $f$ by a suitable complex number, 
it can be written as $f= 1 + \overline{f}$, where $\overline{f}(0,0) \ne 0$.

Suppose that $\pi_p(L_f) \neq 0$. In this case denote by $\M=\{m_1,\ldots,m_k\}$ the set of all monomials in $L_f$ with nonzero coefficients 
such that $\pi_p(m_i) \ne 0$. 

We say that $m_i < m_j$ if there exist $(r,s) \in \Z_+^2 \setminus \{(0,0)\}$ satisfying $x^ry^sm_i = m_j$. Here $\Z_+$ denotes the set of all 
non--negative integers. Observe that $<$ defines a partial order in the finite set $\M$.

Take a minimal element of $\M$ with respect to this partial order and denoted it by $m_*$. Then from the equation
\[
  a\frac{\partial f}{\partial x} + b\frac{\partial f}{\partial y} = L_f + \overline{f}L_f \ ,
\]
one can see that $m_*$ appears with nonzero coefficient on the right hand side, and with zero coefficient on the left hand side
(because $m_* \notin TI_p$). Therefore we obtain a contradiction and the proof is completed.
\cqd

\begin{corollary}\label{ccc}
If $g/f$ is an exponential coefficient of $X$ with $f(p) \ne 0$, then $L_e \in \ker \pi_p$, where $e=\exp(g/f)$.
\end{corollary}

\noindent{\it Proof}: The proof can be made using similar arguments to the previous proposition. The main observation is that 
from the equality
\begin{equation}\label{www}
  X\left(\exp\left(\frac{g}{f}\right)\right) = L_e \exp\left(\frac{g}{f}\right) \ ,
\end{equation}
we obtain
\[
   a\left( f\frac{\partial g}{\partial x} - g\frac{\partial f}{\partial x} \right) + 
   b\left( f\frac{\partial g}{\partial y} - g\frac{\partial f}{\partial y} \right) = L_ef^2= L_e \left(
   f^2(p) + \mbox{h.o.t.} \right) \ .
\]
\cqd

\begin{definition} \rm
We call $\ker \pi_p$ the {\it space of admissible cofactors of invariant algebraic curves of $X$ not passing through $p$}.
\end{definition}

\subsection{The line at infinity}

We say that the {\it line at infinity is invariant} by $X$ if $xb_d - ya_d \ne 0$, where $a_d$ and $b_d$ are the homogeneous
part of degree $d$ of the components of the vector field  (\ref{campo}). This is due to the fact that when we extend the vector
field to a foliation of the complex projective plane $\CP$, the line at infinity turns to be invariant. 

When $xb_d - ya_d \equiv 0$ the vector field $X$ given by (\ref{campo}) can be written in the form
\begin{equation}\label{ttt}
  X = \overline{a} \px + \overline{b} \py + h \left( x \px + y \py \right) \ ,
\end{equation}
where $\overline{a}$ and $\overline{b}$ are polynomials of degree at most $d-1$ and $h$ is a homogeneous polynomial 
of degree $d-1$. In fact $h = a_d /x = b_d /y$. The homogeneous polynomial $h$ defines a set of points in $\mathbb{C}P(1)$, 
this set is formed by the tangencies of the induced foliation with the line at infinity. The relevance of such tangencies for 
our study can be seen in the next proposition.

\begin{proposition} \label{bbb}
Suppose the line at infinity is not invariant. Let $L$ be a cofactor and $L_{d-1}$ be its homogeneous part of highest degree.
Then the following statements hold.
\begin{itemize}
\item[(a)] If $L$ is a cofactor of an invariant algebraic curve $f$ of $X$, then  $L_{d-1} = \deg(f) h$.
\item[(b)] If $L$ is a cofactor of an exponential factor $e=\exp(g/f)$ of $X$, then  $L_{d-1} = 0$.
\end{itemize}
\end{proposition}

\noindent{\it Proof}: Let $k$ be the degree of $f$, and we denote by $f_k$ the homogeneous part of higher degree of
$f$. Since $X(f)=L_f f$ and using Euler's formula for homogeneous polynomials, we get
\[
  \overline{a} \frac{\partial f}{\partial x} + \overline{b}\frac{\partial f}{\partial y}
  + h \left( x\px + y\py \right)(f-f_k) + khf_k = L_f f \ .
\] 
Comparing the homogeneous parts of highest degree in this equality it is easy to see that $khf_k = L_{d-1}f_k$. And consequently 
statement (a) is proved.

Denote by $l$ the degree of $g$, and as usual $g_l$ is the homogeneous part of highest degree of $g$. From (\ref{www}) it follows that 
\[
  ( \overline{a} + hx ) \left( f\frac{\partial g}{\partial x} - g\frac{\partial f}{\partial x} \right) +
  ( \overline{b} + hy ) \left( f\frac{\partial g}{\partial y} - g\frac{\partial f}{\partial y} \right) =
  L f^2 \ .
\]
Equalizing the higher homogeneous parts as before, we are able to write
\[
   h \left( x  \left( f_k\frac{\partial g_l}{\partial x} - g_l\frac{\partial f_k}{\partial x} \right) +
            y  \left( f_k\frac{\partial g_l}{\partial y} - g_l\frac{\partial f_k}{\partial y} \right)   \right) =
  L_{d-1} f_k^2 \ .
\]
Applying Euler's formula again we arrive at the following equation
\[
  (l-k) h f_k g_l = L_{d-1} f_k^2 \ .
\]
Comparing the degrees of both sides we get that the degree of $f$ must be equal to the degree of $g$, i.e., $l=k$. Hence we can conclude
the statement (b).
\cqd
 
\section{Darboux Integrability}\label{darboux}

The next theorem, without taking into account the exponential factors, was proved by Darboux in \cite{D}. This later improvement 
was made in \cite{CL1,CL2}. Its proof follows from straightfoward computations, for more details see \cite{CL2}.

\begin{theorem}\label{classico}
Let $X$ be a vector field. If $X$ admits $p$ distinct invariant algebraic curves $f_i$, for $i = 1, \ldots , p$, and 
$q$ independent exponential factors $e_j$, for $j = 1, \ldots, q$. Then the following statements hold.
\begin{itemize}
\item[(a)] If there are $\lambda_i, \rho_j \in \C$ not all zero such that $\displaystyle{\sum_{i=1}^{p}\lambda_i L_{f_i} +
\sum_{j=1}^{q}\rho_j L_{e_j} = 0 }$, then the (multi--valued) function
$f_1^{\lambda_1}\cdots f_p^{\lambda_p} e_1^{\rho_1}\cdots e_q^{\rho_q}$ is a first integral of the vector field $X$.
\item[(b)] If there are $\lambda_i, \rho_j \in \C$ not all zero such that $\displaystyle{\sum_{i=1}^{p}\lambda_i L_{f_i} +
\sum_{j=1}^{q}\rho_j L_{e_j } = -\mbox{\rm div}  X}$, then the function $f_1^{\lambda_1}\cdots f_p^{\lambda_p} e_1^{\rho_1}\cdots e_q^{\rho_q}$
 is an integrating factor of $X$.
\end{itemize}
\end{theorem}

\begin{definition} \rm
Let $X$ be a vector field of degree $d$, and $\S \subset \C^2$ a finite set of points (eventually empty). The {\it restricted cofactor space 
with respect to $\S$}, $\Sigma_{\S}$, is defined by  
\[
\Sigma_{\S} = \left( \bigcap_{p \in \S} \ker \pi_p \right) \bigcap \Gamma \ ,
\]
where 
\begin{equation*}  
  \Gamma = \begin{cases}
  \C_{d-1}[x,y]               & \mbox{if the line at infinity is invariant}, \\
  \C_{d-2}[x,y] \oplus  h \cdot \C & \mbox{otherwise}.
  \end{cases}
\end{equation*}
We recall that $h$ is the homogeneous polynomial of degree $d-1$ describing the tangencies of the vector field $X$ with the line at  infinity, 
which appears in formula (\ref{ttt}).
\end{definition}

In \cite{CLS} the notion of independent singular points was introduced in order improve the Darboux theory of integrability. 
Recalling that  each point $p=(x_0,y_o)$ in $\C^2$ defines a maximal ideal, $\m_p = \langle x-x_0,y-y_0 \rangle \subset \C[x,y]$.
We can say that $p_1,\ldots ,p_r$  are {\it independent singular points} if 
\[
  \dim \left( \cap_{i=1}^{r} \m_{p_i} \right) \cap \C_{d-1}[x,y] = dim \C_{d-1}[x,y] - r \ .
\]
With this notation, it is easy to prove the following result.

\begin{proposition} \label{melhora}
Assume that $\S=\{p_1,\ldots,p_r\}$ are independent singular points of a vector field of degree $d$. Then the following inequality holds
\[
  \dim \Sigma_{\S} \le \dim \C_{d-1}[x,y] - r \ .
\]
\end{proposition}
 
\begin{theorem}({\bf Integrability Theorem}) \label{main}
Let $X$ be a vector field of degree $d$. Assume that $X$ has $f_1,\ldots, f_q$ distinct invariant algebraic curves, and that $p_1,\ldots, p_r$
are singular points of $X$. Suppose that $f_j(p_k) \ne 0$ for $j=1,\ldots, q$ and $k=1,\ldots, r$. Then following statements hold.
\begin{itemize}
\item[(a)] Let $\mu = \displaystyle{\sum_{j=1}^p \mu_i(X,f_j)}$ and $\sigma = \dim \Sigma_{\S}$, where $\S = \{p_1,\ldots,p_r \}$.
  \begin{itemize}
  \item[(a.1)] If $\mu \ge \sigma +2$, then $X$ has a rational first integral.
  \item[(a.2)] If $\mu \ge \sigma +1$, then $X$ has a Darbouxian first integral.
  \item[(a.3)] If $\mu \ge \sigma$, then $X$ has either a Darbouxian first integral or a Darbouxian integrating factor.
  \end{itemize}
\item[(b)] If there exists a positive integer $l$ such that 
\[ 
   \sum_{j=1}^p \deg (f_j) \mu_{a,l}(X,f_j) \ge n_l(X)  \ ,
\]
then $X$ admits a rational first integral. If we denote the line at infinity as $\L_{\infty}$, then  
\begin{equation*}
 n_l(X) = \begin{cases} 
                      \displaystyle{\frac{d(l^4+6l^3+11l^2+6l)-l^4-2l^3+l^2-6l}{8}}, & \mbox{if $\L_{\infty}$ is invariant}, \\  & \\
                      \displaystyle{\frac{d(l^4+6l^3+11l^2+6l)-2l^4-8l^3-10l^2-4l}{8}}, & \mbox{otherwise}.
        \end{cases}
\end{equation*}           
\end{itemize}
\end{theorem}

\noindent{\it Proof}: The hardest part of the proof has already been made in the previous sections, and here this is 
reflected in the way that the assumptions in the statements of the theorem are presented. Now to prove parts (a.2) and (a.3)
is just a matter of counting dimensions and using Theorem \ref{classico}. More precisely, one has just to observe that
all possible cofactors are contained in $\Sigma_{\S}$ (see Propositons \ref{aaa}, \ref{bbb} and Corollary \ref{ccc}). Hence if $\mu$ is at least
$\sigma$ we can assure that the hypothesis of statement either (a) or (b) of Theorem \ref{classico} are fullfilled, and hence
(a.3) follows. 

In the case that $\mu$ is at least 
$\sigma +1$ we can garantee that the assumption of statement (a) of Theorem \ref{classico} holds, therefore (a.2) holds.

When $\mu$ is at least $\sigma + 2$, from (a.2) we can obtain two independent Darbouxian first integrals, say $H_1$ and $H_2$. We can see easily
that the integrating factor $R_i$ associated to $\log H_i$ is a rational function. Since the quotient of two integrating factors is a first integral,
the statement (a.1) follows from the independence of $H_1$ and $H_2$.

Statement (b) can be proved observing that the degree of $\E_l(X)$, together with Theorem \ref{ddd} and Proposition \ref{pext}
give bounds for the number of invariant algebraic curves of degree at most $l$. For more details see \cite{P}.
\cqd

\begin{corollary}
If there exists a positive integer $l$ such that 
\[ 
   \sum_{j=1}^p \deg (f_j) \mu_{g,l}(X,f_j) \ge n_l(X)  \ ,
\]
then $X$ admits a rational first integral.
\end{corollary}

\section{Strong Algebraic Multiplicity}\label{strong}

\begin{definition}\rm
If $X$ is a vector field on $\C^2$ the {\it $n$--th extactic ideal of $X$}, $\E I_n(X)$, 
is generated by  $\sigma_{(k_1,\ldots,k_l)}$, where   
\begin{equation}
\mbox{det} \begin{pmatrix}  X^{k_1}(v_1)    &    X^{k_1}(v_2)       & \cdots & X^{k_1}(v_l)          \cr
                           X^{k_2}(v_1)       & X^{k_2}(v_2)       & \cdots & X^{k_2}(v_l)       \cr
                           \vdots       & \vdots       & \cdots & \vdots       \cr
                           X^{k_l}(v_1) & X^{k_l}(v_2) & \cdots & X^{k_l}(v_l) \end{pmatrix},
\end{equation}
where  $0 \le k_1 < k_2 < \cdots < k_l$ and the $k_i$'s are integers,  $v_1, v_2, \ldots, v_l $ 
is a basis of  $\C_n[x,y]$  (so $l=(n+1)(n+2)/2$).
\end{definition}

\begin{definition}\rm
An invariant algebraic curve $f$ of degree $n$ for the vector field $X$ has {\it $l$--strong algebraic multiplicity $m$} when 
$m$ is the smallest positive integer such that the
$m$--th power of $f$ belongs to the extactic ideal, $\E I_l(X)$. We denote $m$ by $\mu_{sa,l}(X,f)$,
or if the context is clear $\mu_{sa,l}(f)$.
\end{definition}

Proposition \ref{ineq1} can be easily generalized to give the pair of inequalities,
\begin{equation}\label{iii}
\mu_{g,l}(X,f) \leq \mu_{sa,l}(X,f) \leq \mu_{a,l}(X,f) \ .
\end{equation}

\begin{example}\label{decide} \rm
We have showed, in Example \ref{dana1alg}, that the square of the invariant straight line $1 - by$ of the vector field $X_{(0,b,d)}$
has $2$--algebraic multiplicity $3$, and hence we know that its $2$--geometric multiplicty is either $2$ or $3$. 
With just two generators of the
second extactic ideal $\E I_n(X_{(0,b,d)})$, namely $\sigma_{(0,1,\ldots,5,6)}$ and $\sigma_{(0,1,\ldots,5,7)}$, one can see that the 
$2$--strong algebraic multiplicity is exactly $2$. Here we use the inequality (\ref{iii}), and also from it follows that the $2$--geometric
multiplicity is exactly $2$.
\end{example}

\section{Final Remarks} \label{final}

Although we studied vector fields in the complex plane, all results here extend easily to real polynomial vector
fields defined on the real plane. To do that one has just to think the real vector fields as been complex, and after 
obtaining a complex first integral or complex integrating we can return to the real world using the results presented
in Section $3$ of \cite{Mac}.

We also remark that it is possible to generalize our results to codimension $1$ foliations on $\C^n$ given by polynomial
$1$--forms. To see a hint how to do that we suggest that reader consult the work of Jouanolou \cite{J}.

To conclude we would like to emphasize that at first sight the concept of geometric multiplicity is not computable, neverless
the algebraic and strong algebraic multiplicity gives a computational approach to bound such multiplicity. We believe that in some
sense the strong algebraic multiplicity contains all information necessary to obtain the geometric multiplicity.

\section{Acknowledgments}

We want to thank the organizers of the {\it Workshop on asymptotic series, differential algebra and finiteness problems
in non--linear dynamical systems} held at the Centre de Recherces Math\'ematiques -- Universit\'e de Montr\'eal. Their support
during our stay allowed us to have a great atmosphere to conceive this paper. Specially we want to express our
gratitude to Dana Schlomiuk and Nicolas Vulpe for providing examples that permitted us to test the different notions of multiplicity.

The first author is partially supported by DGES grant number PB96--1153 and by CICYT grant number 1999SGR
00399.

\end{document}